\magnification=1000

\input local.layout
\input macros/math.macros

\proofmodefalse

\def\Dots{s^{\k{-3.5}\raise2.2pt\hbox{$\scriptscriptstyle \circ$}}}
\def\Dotw{w^{\k{-5}\raise2.2pt\hbox{$\scriptscriptstyle \circ$}}}

\def\Ddotl{1^{\k{-5.5}\raise2.2pt\hbox{$\scriptscriptstyle \triangle$}}}
\def\Ddots{s^{\k{-5}\raise2.2pt\hbox{$\scriptscriptstyle \triangle$}}}
\def\Ddotw{w^{\k{-6.5}\raise2.2pt\hbox{$\scriptscriptstyle \triangle$}}}
\def\Ddotx{x^{\k{-5}\raise2.2pt\hbox{$\scriptscriptstyle \triangle$}}}
\def\Ddotnull{\null^{\k{-1}\raise2.2pt\hbox{$\scriptscriptstyle \triangle$}}}

\def\dots{{s^{\k{-3.5}\raise2.2pt\hbox{$\scriptscriptstyle \bullet$}}}}
\def\dotw{{w^{\k{-5}\raise2.2pt\hbox{$\scriptscriptstyle \bullet$}}}}
\def\dotx{{x^{\k{-4}\raise2.2pt\hbox{$\scriptscriptstyle \bullet$}}}}
\def\doty{{y^{\k{-3.5}\raise2.2pt\hbox{$\scriptscriptstyle \bullet$}}}}
\def\dott{{\null^{\k{-3.5}\raise2.2pt\hbox{$\scriptscriptstyle \bullet$}}}}

\def\ph{\phantom{-}}

\def\term#1{\dbll #1 \dblr}
\def\x{\times}
\def\Vee{{\scriptscriptstyle \vee}}

\def\frac#1#2{#1/#2}

\def\magma/{{\tt MAGMA}}
\def\seck/{\S 3}
\def\secc/{\S 1}
\def\sect/{\S 2}
\def\seccomp/{\S 4}

\def\e{{\goth k}}
\def\f{{\goth e}}

\def\dbll{\langle\k{-1.8}\langle}
\def\dblr{\rangle\k{-1.8}\rangle}
\def\Ss{S}
\def\ad{\k{1}\raise1pt\hbox{$\scriptscriptstyle \diamondsuit$}\k{0.5}}
\def\ext{{\rm ext}}
\def\vare{\varepsilon}
\def\W{W}
\def\height{{\rm ht}}


{\bigbold A simple way to compute structure constants of semi-simple Lie algebras}

\bigskip

Bill Casselman\hfill\break
University of British Columbia\hfill\break
{\tt cass@math.ubc.ca}

\bigskip

Suppose ${\goth g}$ to be a complex Lie algebra
with basis $(x_{i})$.  Then
$$ [x_{i}, x_{j}] = {\sum} _{k} c^{k}_{i,j} x_{k} $$
for some numbers $c^{k}_{i,j}$, called {\sl structure constants}.

For the classical semi-simple groups, they can be computed using the defining
matrix representations, but for exceptional groups something else is needed.
This is by no means a trivial exercise. 
One commonly used technique for computing 
them, which has been used in the computer program \magma/, is described 
well by \refto{Cohen-Murray-Taylor}{2005}.  It 
relies on the additive structure of roots. Another
method, that works only for simply-laced root systems and relies on associated
affine root systems, is explained in \refto{Frenkel-Kac}{1980}.
A version of this for the remaining root 
systems can be found in \refto{Rylands}{2000}.  
It uses the identification of these others as folded 
quotients of simply laced ones.

In \refto{Casselman}{2015b} I explained yet another way to compute these structure
constants by implementing an idea
originally found in \refto{Tits}{1966a}.  Tits' idea was to
replace the additive structure by features of the normalizer of a maximal torus.
This introduced some mathematical structure 
to the problem of computing structure constants 
that was missing in the standard approach.
In practice, computation based on this method
went fairly rapidly and seemed at least roughly comparable in efficiency to 
reported runs of the standard computation. 
There were, however, a number of rather ugly and presumably inefficient
formulas involved in this new algorithm.  In May of 2014,
it was suggested by Robert Kottwitz, that an observation of his
about choosing bases of semi-simple 
Lie algebras might make it possible to bypass the
nastiest parts in a more elegant manner.
In this paper, with Kottwitz' permission, I'll explain how this goes. 

\ignore {As I say, the principal ideas in this paper are due to Kottwitz.  
The remarks on practical implementation, including
details on how to use Tits' results for computation, are due to me.}

Kottwitz' basic observation is very simple, and can be briefly summarized. Suppose
$$ \eqalign {
	G \k{1} &= \hbox{a simple, connected, simply connected, complex group}\cr
	{\goth g}\k{2}  &= \hbox{Lie algebra of $G$} \cr
	B \k{1}  &= \hbox{Borel subgroup}\cr
	T \k{1}  &= \hbox{maximal torus in $B$}\cr
	\Sigma \k{1}  &= \hbox{associated root system}\cr
	\Delta &= \hbox{associated simple roots}\cr
	W \k{-2} &= \hbox{Weyl group.}\cr
} $$
Because $G$ is simply connected, the coroot lattice
$X_{*}(T)$ may be identified with the lattice
spanned by the simple coroots $\alpha^{\Vee}$.

The root spaces ${\goth g}_{\gamma}$ all have dimension one.
Fix for each $\alpha$ in $\Delta$ a spanning element
$e_{\alpha}$ in ${\goth g}_{\alpha}$.  The triple $(B, T, \{ e_{\alpha}\})$ make
up a \defn{frame} for $G$.  The set of all frames
is a principal homogeneous space for the 
adjoint quotient of $G$. (This notion originated in
work of French mathematicians.  In French the term is
`\'epinglage', which some translate literally into the noun `pinning'.  
But `frame' is the term adopted in the English translation of 
Bourbaki's treatise on Lie algebras.)  The point is that computation in
the group or its Lie algebra must start with
a frame, which amounts to a kind of basis of the Lie algebra.

As I'll recall later, Chevalley has defined integral structures on ${\goth g}$
and $G$.  The map 
$$ \{ \pm 1 \}^{\Delta} \longrightarrow T, 
	\quad (c_{\alpha}) \longmapsto \prod_{\alpha \in \Delta} \, \alpha^{\Vee}(c_{\alpha})$$
identifies $T({\Bbb Z})$ with a two-torsion group.
If ${\cal N}({\Bbb Z})$ is the group
of integral points in the normalizer ${\cal N} = N_{G}(T)$,
 it fits into a well understood extension
$$ 1 \longrightarrow T({\Bbb Z})
\longrightarrow {\cal N}({\Bbb Z})
\longrightarrow W
\longrightarrow 1 \, . $$
\refto{Tits}{1966b} associates to a given frame a certain 
convenient section $w \mapsto \dotw$\k{4} of the last quotient map,
and described this extension precisely
enough to enable computations in it.
This extension certainly does not generally split
(as it does for $\GL_{n}$, say).
But now let $V_{\Bbb Z}$ be the direct sum of non-trivial root
spaces in ${\goth g}_{\Bbb Z}$.   
Let $\Ss({\Bbb Z})$ be the subgroup of transformations in
$\GL(V_{\Bbb Z})$ that act as $\pm 1$ on each root space.
It may be identified with $\Hom(\Sigma, \pm 1) = (-1)^{\Sigma}$.
The adjoint action of $T$
defines a canonical homomorphism from $T({\Bbb Z})$ 
to $\Ss({\Bbb Z})$: $\alpha^{\Vee}(x)$
goes to $(x^{\<\gamma, \alpha^{\Vee}>})_{\gamma\in\Sigma}$.  
The kernel is $Z_{G}$.
The homomorphism from $T$ to $\Ss$ gives rise to an extension
$$ 1 \longrightarrow \Ss({\Bbb Z}) = \{ \pm 1 \}^{\Sigma}
\longrightarrow {\cal N}_{\rm ext} ({\Bbb Z})
\longrightarrow W
\longrightarrow 1 \, . $$
Although it does not act as automorphisms of ${\goth g}$,
the extension does act on $V_{\Bbb Z}$, compatibly with the adjoint 
action of $T({\Bbb Z})$.   
Kottwitz' notable observation is that {\sl this new extension splits}, 
and he gives an explicit splitting $w \mapsto \Ddotw$.
It has the property that if $w\lambda = \lambda$ then $\Ddotw$ 
acts as the identity on ${\goth g}_{\lambda}$.  This allows
one to specify a natural choice of Chevalley basis, once
one fixes one $e_{\gamma}$ in each $W$-orbit in $\Sigma$.
One consequence of the new method is a very simple description of
the action of ${\cal N}({\Bbb Z})$ on ${\goth g}$.  This is especially important
in applications to computation in the group $G$ rather than just its Lie algebra.

Some of the previous methods 
known have the virtue that they may be extended to all
Kac-Moody root systems (see \refto{Casselman}{2015b}).
Some variant of the method I describe here will work for
a large class of these.  I do not see how it can be extended to all of them, but 
one might hope that some variation of Kottwitz' idea will work, taking into
account some explicit obstruction.
One promising prerequisite for 
extending the method to Kac-Moody algebras can be found 
in \refto{Carbone \etal/}{2015}, which classifies 
conjugacy classes of simple roots.
 
Curiously, it was in \refto{Langlands-Shel\-stad}{1987}
that an explicit formula for a defining $2$-cocycle
of Tits' sections $w \mapsto \dotw$ first appeared.
Recently, Tasho Kaletha has found other 
applications of Tits' construction and results of this paper
to related problems.
I wish to thank him  for comments on an earlier version.

\contents

For $g$ in $G$, $x$ in ${\goth g}$, I'll write
$$ g \ad x = \Ad(g) x \, . $$
Even though ${\cal N}_{\ext}(\Bbb Z)$
does not act by automorphisms of ${\goth g}$,
I'll use this notation for its action on $V_{\Bbb Z}$ as well.

I'll usually refer to \refto{Tits}{1966a} as [T].

\def\z{\zero}

\section{Chevalley bases}
Fix once and for all a maximal torus $T$ in $G$.
The associated roots are the non-trivial
characters by which $T$ acts on eigenspaces,
each of which has dimension one.  For the moment,
suppose $\gamma$ to be a root.  If $e \ne 0$
lies in ${\goth g}_{\gamma}$, then for every $f$ in
${\goth g}_{-\gamma}$ the bracket $h = [f, e]$
will lie in the Lie algebra ${\goth t}$.  There will exist exactly one 
$f$ such that $[h,e] = 2 e$.
In these circumstances I'll call $(e, h, f)$
an \defn{$\SL_{2}$ triple}.
It is completely determined by the choices of $T$ and of $e$
in ${\goth g}_{\gamma}$.

Given such a triple, there exists a unique embedding
$\iota_{e}$ of $\SL_{2}$ into $G$ whose differential takes
$$ \eqalign {
	\mmatrix { \k{2} \z & \k{2} 1 \cr
			\k{2} \z & \k{2} \z \cr } &\longmapsto e \cr
	\mmatrix { \k{-2} \phantom{-}\z & \k{-2} \z \cr
			 \k{-2} -1 & \k{-2} \z } &\longmapsto f \cr
	\mmatrix { 1 & \k{-4} \phantom{-}\z \cr
			 \z & \k{-4} -1 } &\longmapsto h \, . \cr
} $$
If we change $e$ to $xe$ with $x \ne 0$, then
$f$ changes to $x^{-1}f$, and 
$\iota_{e}$ changes to its conjugate by 
$$ \mmatrix{ \sqrt{x} & \k{-5} \z \cr \z & \k{-5} 1/\k{-2}\sqrt{x} \cr } \, . $$
The associated embedding of ${\Bbb C}^\x$ is 
the coroot $\gamma^{\Vee}$, and is independent of the choice of $e$. 

Now fix in addition a Borel subgroup $B$ containing $T$.
Let $\Delta$ be the corresponding set of simple roots 
and for each $\alpha$ in $\Delta$ fix an element $e_{\alpha} \ne 0$
in ${\goth g}_{\alpha}$.
The triple $(B, T, \{ e_{\alpha} \})$ makes up a frame
for $G$.  The set of all frames is a principal homogeneous
space for the group of inner automorphisms of $G$.

The frame determines embeddings $\iota_{\alpha}$
of $\SL_{2}$ into $G$, one for each
simple root.  Let $h_{\alpha}$ be the image under $d\iota_{\alpha}$ 
of 
$$ \mmatrix { 1 & \k{-3}\phantom{-}\zero \cr \zero & \k{-3}-1 \cr } \, . $$
The image $e_{-\alpha}$ of 
$$ \mmatrix { \k{-2}\phantom{-}\zero & \k{1}\zero \cr
		\k{-2}-1 & \k{1}\zero \cr }
$$
is the unique element of ${\goth g}_{-\alpha}$ such that
$$ [e_{-\alpha}, e_{\alpha} ] = h_{\alpha} \, . $$
This choice of sign is Tits'. 
It is not the common one, 
but it is exactly what is needed to make
his analysis of structure constants work.
The point is that there exists an 
automorphism $\theta$ of ${\goth g}$
acting as $-I$ on ${\goth t}$
and taking each $e_{\alpha}$ ($\alpha \in \Delta$) to $e_{-\alpha}$. 
It is uniquely determined by the choice of frame.

\subsec{Tits' section}
The group ${\cal N}({\Bbb Z})$ fits into a short exact sequence 
$$ 1 \longrightarrow T({\Bbb Z}) = (\pm 1)^{\Delta} 
\longrightarrow {\cal N}({\Bbb Z})
\longrightarrow W
\longrightarrow 1 \, , $$
and \refto{Tits}{1966b} shows 
how to define a particularly convenient section.
Define 
$$ \dots_{\alpha} = \iota_{\alpha}\left( \mmatrix { \phantom{-}\zero & 1 \cr
			-1 & \zero \cr }\right) \, . $$
It lies in the normalizer of $T$.
Suppose $w$ in $W$ to have the reduced expression $w = s_{1} \ldots s_{n}$.
Then the product 
$$ \dotw = \dots_{1} \ldots \dots_{n} $$
depends only on $w$, not the particular product 
expression. The defining relations
for this group, given those for $T$ and $W$, are
$$ \eqalign {
	(xy)\k{2}\dott &= \dotx \doty \quad (\hbox{when } \ell(xy) = \ell(x) + \ell(y)) \cr
	\dots^{\k{1}2}_{\alpha} &= \alpha^{\Vee}(-1) \quad (\alpha \in \Delta) \, . \cr
} $$
The effect of a change of basis is easy to figure out.

\lemma{change}{Suppose $f_{\alpha} = c_{\alpha}e_{\alpha}$, and 
let $\Dots_{\alpha}$ be the corresponding element of ${\cal N}({\Bbb Z})$.
Then
$$ \Dots_{\alpha} \ad x_{\lambda} 
	= (c_{\alpha})^{\<\lambda, \alpha^{\Vee}>} \dots_{\alpha} \ad x_{\lambda} \, . $$
}

\subsec{Chevalley's formula} 
Suppose $(e_{\gamma}, h_{\gamma}, e_{-\gamma})$
to form an $\SL_{2}$ triple, and suppose
that $e_{\gamma}^{\theta} = c e_{-\gamma}$.  If 
$f_{\gamma} = e_{\gamma}/\sqrt{c}$ and $f_{-\gamma} = \sqrt{c} e_{\gamma}$,
then $(f_{\gamma}, h_{\gamma}, f_{-\gamma})$ also make up
an $\SL_{2}$ triple,
with $f_{\gamma}^{\theta} = f_{-\gamma}$.  Up to sign---but only up
to sign---$f_{\gamma}$ is unique with this invariance condition.

Any complete set $\{ e_{\gamma} \}$ invariant under $\theta$ up to sign
is often called a Chevalley basis (with respect to the given frame).
It determines an integral structure on the Lie algebra ${\goth g}$.

\definition{integral-basis}{I'll call such a 
basis an \defn {integral basis}.
If it is actually invariant under $\theta$, as it is here,
I'll call it an \defn{invariant} basis.}

{\bold Remark.} The integral structure on ${\goth g}$
is determined by the frame, and more directly from the involution $\theta$
it defines.  It is curious that $\theta$ also determines a maximal
compact subgroup of $G$.  Of course for $p$-adic groups,
there is a more immediate relation between integral
structure and compact subgroups.

\endremarks

Given any integral basis, Chevalley proved that if $\lambda$, $\mu$, $\nu$
are roots with $\lambda + \mu + \nu = 0$ then 
$$ [e_{\lambda}, e_{\mu}] = \pm (p_{\lambda,\mu} + 1) e_{-\nu} \,. \equation{cheveq} $$
Here $p_{\lambda,\mu}$ is the least $p$ such that $\mu - p\lambda$ is
a root.  This was the crucial result used to construct the Chevalley
groups over arbitrary fields.

The possible values for the string constants $p_{\lambda,\mu}$
(associated to finite root systems) are shown in the following figures:

\figures{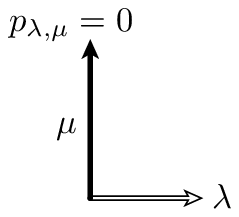}{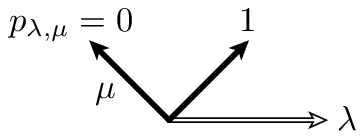}

\figures{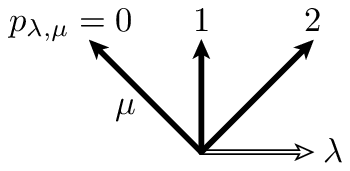}{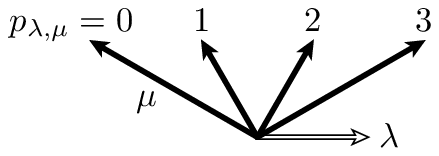}

The fourth figure
occurs only in type $G_{2}$.  In practice, we shall be interested in computing
$p_{\lambda,\mu}$ only when $\< \mu, \lambda^{\Vee}> \le 0$.  Under this assumption,
as the figures illustrate:
$$ p_{\lambda,\mu} = \cases {
	0 & if $\mu - \lambda$ is not a root \cr
	1 & otherwise \cr
} \quad (\hbox{assuming } \< \mu, \lambda^{\Vee}> \le 0 ) \, .  \equation{plambdamu} $$

I refer to \refto{Chevalley}{1955} or \refto{Carter}{1972} for the original 
proof of \lreftosatz{cheveq} and to \refto{Casselman}{2015a} 
for a proof extracted from [T],
which works uniformly for all Kac-Moody groups.  Tits'
choice of the $e_{-\alpha}$ (as opposed to the more common choice
with the opposite sign) introduces an elegant symmetry
that greatly simplifies both proofs and formulas.

{\bold Remark.} Ultimately, Chevalley's formula depends on the simple
fact that for strings of length $2$, as in the second figure above, 
one always has $\Vert \lambda \Vert \ge \Vert \mu \Vert$.  
That is to say, the following configuration never occurs.

\figure{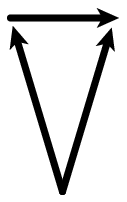}

\endremarks

Determining the sign in \lreftosatz{cheveq} has always seemed
rather mysterious.  Of course there can be no
simple formula, since the choice of an integral basis
is not canonical.  But I don't think it has ever been
very clear what is going on.  Changing
even one $e_{\gamma}$ to $-e_{\gamma}$ forces a lot of other sign changes
without apparent pattern.
The situation has now been cleared up somewhat by Kottwitz,
who has explained to me how to choose an almost canonical
integral basis.  I'll discuss that in \seck/.

In \seccomp/ I'll show how Kottwitz' basis simplifies the computation
of the signs in Chevalley's formula.  The starting
point, at least, is the same as it was in \refto{Casselman}{2015b},
in which I have already outlined the principal ingredients of a recipe
for the constants.  One of the troublesome
points in the earlier approach was a somewhat arbitrary
choice of integral basis.  Kottwitz' basis
eliminates this inconvenience.

\ignore{I have fixed the frame $\{e_\alpha\}$, but I have not extended it to
an integral basis.  Curiously, Kottwitz' method specifies an 
invariant integral basis that is related 
to this original frame but does not agree with it.}

\section{Tits' idea}
In order to understand how Kottwitz' basis makes calculation
of structure constants simple, I must explain how
it fits into the scheme covered in \refto{Casselman}{2015b} for computing
structure constants.  I'll do that in this section and the next.

In this one I shall recall results of Tits alluded to at the beginning of \S 1.
We have seen there that a choice of root vector $e$ 
determines an embedding $\iota_{e}$ of $\SL_{2}$ into $G$, and in particular
determines the element 
$$ \sigma_{e} = \iota_{e} \left( \mmatrix { \phantom{-}\z & 1 \cr
			-1& \z \cr } \right) \, . $$
Tits starts with a variation on this fact, an elementary
observation about $G = \SL_{2}$. Let $T$ be the subgroup of
diagonal matrices.  Its normalizer in $G$ is the union of $T$
itself and the subset $M$ of matrices of the form
$$ \mmatrix { 
\k{-2} \phantom{-}\z & \k{-3} x \cr
\k{-2} -1/x & \k{-3} \z } \, . $$
Let ${\goth g}_{+}$ be the Lie algebra of upper nilpotent
matrices
$$ \mmatrix { \z & x \cr
			\z & \z } \, , $$
${\goth g}_{-}$ that of lower nilpotent ones.
The following is Proposition 1 of \S 1.1 of \refto{Tits}{1966a}:

\lemma{tits-cts}{Suppose $e_{+}$ in ${\goth g}_{+}$,
$e_{-}$ in ${\goth g}_{-}$, $\sigma$ in $M$.  The following are equivalent:
$$ \eqalign {
	(a) \hbox{\hskip3em} & \exp(e_{+}) \exp(e_{-}) \exp(e_{+}) = \sigma \, ; \cr
    (b) \hbox{\hskip3em} & \exp(e_{-}) \exp(e_{+}) \exp(e_{-}) = \sigma \, . \cr
} $$
If any one of these three matrices is specified,
conditions (a) or (b) determine the other two uniquely.} 

\proof/. An easy matrix calculation shows that if
$$ \mmatrix { 1 & x \cr
			\z & 1 \cr }
	\mmatrix { 1 & \z \cr
			y & 1 \cr }
	\mmatrix { 1 & x \cr
			\z & 1 \cr }
$$
lies in the normalizer of $T$, then $y = -1/x$, in which case the product is
$$ \mmatrix { \phantom{-}\z & x \cr
			-1/x & \z } \, . $$
This proves the last claim. The equivalence of (a) and (b)
follows from the equation
$$ 
	\mmatrix { \k{-2}\phantom{-}\z & \k{-3}x \cr
			\k{-2}-1/x & \k{-3}\z } 
	\mmatrix { 1 & x \cr
			\z & 1 \cr }
	\mmatrix { \z & \k{-6} -x \cr
			1/x & \k{-6} \phantom{-}\z } 
	= 
	\mmatrix { \k{-2}\phantom{-}1 & \z \cr
			\k{-2}-1/x & 1 \cr } \, . \mendproof $$
I'll call the triplet $(e_{+}, \sigma, e_{-})$ \defn{compatible}, and sometimes
express the element in the normalizer as $\sigma_{e_{+}}$,
which is the same as $\sigma_{e_{-}}$.
This is one place where Tits' choice of $e_{-}$, rather then $-e_{-}$,
is significant.

Now let $G$, $T$ be arbitrary, as earlier.
Suppose given some 
$f_{\gamma}$ generating ${\goth g}_{\gamma} \cap {\goth g}_{\Bbb Z}$.
It is unique up to sign.
As we have seen, it determines
an embedding of $\SL_{2}({\Bbb Z})$ into $G({\Bbb Z})$ and elements 
$h_{\gamma}$, $f_{-\gamma}$ spanning a copy of ${\goth sl}_{2}$.
We also get then an element $\sigma_{\gamma}$ in ${\cal N}({\Bbb Z})$,
the image of
$$ \mmatrix { \phantom{-}\zero & 1 \cr
			-1 & \zero \cr } \, . $$
These also satisfy the equation
$$ \exp(f_{\gamma}) \exp(f_{-\gamma}) \exp(f_{\gamma}) = \sigma_{\gamma} \, . $$
I'll also call the triplet $(f_{\gamma}, \sigma_{\gamma}, f_{-\gamma})$
compatible.   If $\gamma = \alpha$
and $f_{\alpha} = e_{\alpha}$, then $\sigma_{\alpha} = \dots_{\alpha}$,
but I do not assume this to hold.
In any case the image of $\sigma_{\gamma}$ 
in $W$ will be $s_{\gamma}$, the reflection 
corresponding to $\gamma$.  The basic observation of Tits ([T], Proposition 1)
is that each of the objects $f_{\pm \lambda}$, $\sigma_{\lambda}$ determines
the other two.  The choice of sign for any one of these
determines a change of sign in the others.

In other words, the choice of an invariant
basis is equivalent to a certain choice of elements
in the normalizer ${\cal N}({\Bbb Z}) = N_{G}(T) \cap G({\Bbb Z})$.

An elementary calculation in $\SL_{2}$ tells us that
$$ \exp(\vare f_{\mu}) \exp(\vare f_{-\mu}) \exp(\vare f_{\mu}) 
	= \mu^{\Vee}(\vare) \sigma_{\mu} \, . \equation{fef} $$
Hence:

\lemma{compatibility}{If $(e_{\lambda}, \sigma_{\lambda}, e_{-\lambda})$
are compatible, then so are
$(c e_{\lambda}, \lambda^{\Vee}(c) \sigma_{\lambda}, ce_{-\lambda})$.}

{

\parindent=2em
\narrower

\noindent
{\sl For the indefinite future, fix an invariant Chevalley basis $(f_{\gamma})$.}

}

I repeat that 
I do {\sl not} assume that $f_{\alpha} = e_{\alpha}$ for simple roots $\alpha$.
This determines also for each $\gamma$ an element $\sigma_{\gamma}$,
subject to the equations
$$ \eqalign { 
	 \sigma_{\gamma}^{-1} &= \gamma^{\Vee}(-1)\sigma_{\gamma} \cr
	\sigma_{-\gamma} &= \sigma_{\gamma} \, . \cr
} $$

Let $M_{\gamma}({\Bbb Z})$ be the subset of $N_{G}({\Bbb Z})$
in the image of $M \subset \SL_{2}$ 
determined by $\gamma$.  It has two elements,
and contains precisely the $\gamma^{\Vee}(\pm 1) \, \sigma_{\gamma}$.

One practical consequence of Tits' observation is this:

\satz{practical}{Lemma}{Suppose $\omega$ to be in ${\cal N}({\Bbb Z})$.
Let $w$ be its image in $W$, and assume that $w\lambda = \mu$.  
$$ \omega \ad f_{\lambda} = \varepsilon \, f_{\mu} $$
if and only if
$$ \omega \k{1} \sigma_{\lambda} \omega^{-1} = \mu^{\Vee}(\varepsilon) \, \sigma_{\mu} \, . $$
}

Here $\varepsilon$ is necessarily $\pm 1$.

\proof/. Because of \lreftosatz{fef}.\endproof

I remind you that the problem we are considering is this:

{

\parindent=2em
\narrower

\noindent
{\sl Given the integral basis $(f_{\gamma})$, 
we want to figure out how to calculate the sign in Chevalley's formula}
$$ [f_{\lambda},f_{\mu}] = \pm (p_{\lambda,\mu}+1) f_{\lambda + \mu} \, . $$

}

Tits has introduced a convenient symmetry into this problem
by his choice of $f_{-\gamma}$.  For example, 
since this basis is invariant under $\theta$, the constants are now the same
for $-\lambda$, $-\mu$ and $\lambda$, $\mu$.
Tits has introduced a second symmetry 
by another simple notion.  I define
a \defn{Tits triple} to be a set of roots
$\lambda$, $\mu$, $\nu$ whose sum is $0$.
He makes this choice instead of taking, 
more conventionally, $\lambda + \mu = \nu$.

In any finite irreducible root system
there are at most two lengths.  Hence if $\lambda + \mu + \nu = 0$,
two of them must be of the same length.
As I have already mentioned, the common length
cannot be greater than the third.  Therefore any Tits triple
can be cyclically permuted to satisfy the condition
$$ \Vert \lambda \Vert \ge \Vert \mu \Vert = \Vert \nu \Vert \, . $$
In this case, I shall call it an \defn{ordered} triple.

\satz{tits-tfae}{Proposition}{{\rm ([T], Lemme 1 of \S 2.5)} Suppose
$(\lambda, \mu, \nu)$ to be a Tits triple.
The following are equivalent:
\medbreak

\beginitems{3em}

\itema it is an ordered triple;

\itemb $s_{\lambda}\mu = -\nu$;

\itemc $\< \mu, \lambda^{\Vee}> = -1$.

\enditems

}

The upshot of the discussion so far is that there  exists a function 
$\varepsilon(\sigma_{\lambda}, \sigma_{\mu}, \sigma_{\nu})$, 
defined on all products 
$M_{\lambda}({\Bbb Z}) \times M_{\mu}({\Bbb Z}) \times M_{\nu}({\Bbb Z})$
whenever $(\lambda, \mu, \nu)$ is a Tits triple, such that
$$ [f_{\lambda}, f_{\mu}] 
	= \varepsilon(\sigma_{\lambda}, \sigma_{\mu}, \sigma_{\nu}) \, (p_{\lambda,\mu} +1)f_{-\nu} \, . $$
Of course I am assuming that the $\sigma$ and $f$ are compatible.
The following is the basis of computation of structure
constants by Tits' method.  

\satz{tits-f}{Proposition}{{\rm ([T], \S 2.9)}  The function 
$\varepsilon(\sigma_{\lambda}, \sigma_{\mu}, \sigma_{\nu})$ 
satisfies these basic properties
\smallbreak

\beginitems{3em}

\item{($\varepsilon$a)} replacing $\sigma_{\lambda}$ by $\sigma_{\lambda}^{-1}$
changes the sign;

\item{($\varepsilon$b)} it is skew-symmetric in any pair;

\item{($\varepsilon$c)}\k{1} it is invariant under cyclic rotation of the arguments;

\item{($\varepsilon$d)} if $\lambda$, $\mu$, $\nu$ are an ordered triple 
with $\sigma_{\lambda}\sigma_{\mu}\sigma_{\lambda}^{-1} = \sigma_{\nu}$ then
$$ \varepsilon(\sigma_{\lambda}, \sigma_{\mu}, \sigma_{\nu})
	= (-1)^{p_{\lambda,\mu}} \, . $$

\enditems

}

The first two are immediate, but the third is not quite so.
Together, these mean that we can apply a permutation to
any triple to reduce to a special case, but
what is now needed is one explicit formula
in that special case---\ie/ to pin down signs. That is what the last does.
It follows from an analysis (in [T], \S 1.3) of the action of
copies of $\SL_{2}({\Bbb Z})$ on the spaces in ${\goth g}$
determined by root strings in ${\goth g}$.

For an ordered triple, because of (b) and the equality of
$\sigma_{\gamma}$ and $\sigma_{-\gamma}$:
$$ \sigma_{\lambda}\sigma_{\mu}\sigma_{\lambda}^{-1} 
	= \nu^{\Vee}(\pm 1) \sigma_{\nu} \, . $$

\satz{rotation}{Proposition}{Suppose
$(\lambda, \mu, \nu)$ to be an ordered Tits triple,
$\varepsilon = \pm 1$.
The following are equivalent:
\smallbreak

\beginitems{3em}

\itema $ \sigma_{\lambda}\sigma_{\mu} \sigma_{\lambda}^{-1} 
	= \nu^{\Vee}(\varepsilon) \, \sigma_{\nu}$;

\itemb $\sigma_{\lambda} \ad f_{\mu} = \varepsilon \, f_{-\nu}$;

\itemc $ [f_{\lambda}, f_{\mu}] 
	= \varepsilon \, (-1)^{p_{\lambda,\mu}} \, (p_{\lambda,\mu}+1) \, f_{-\nu}$.

\enditems
}

Combining these two propositions:

\theorem{abc-cor}{Suppose $(\lambda, \mu, \nu)$ to be an ordered triple,
$\varepsilon = \pm 1$.  Assume that 
$$ \sigma_{\lambda} \ad \f_{\mu} = \varepsilon \k{1} f_{-\nu} \, . $$
If 
$$ c =  \varepsilon \k{1} (-1)^{p_{\lambda,\mu}} $$
then 
$$ \eqalign {
	[f_{\lambda}, f_{\mu}] &=  c \k{1} (p_{\lambda,\mu} + 1)f_{-\nu} \cr
	[f_{\mu}, f_{\nu}] &= c \k{1} (p_{\mu,\nu} + 1)f_{-\lambda} \cr
	[f_{\nu}, f_{\lambda}] &= c \k{1} (p_{\nu,\lambda} + 1)f_{-\mu} \, . \cr
} $$
}

This will be the basis of computations, once we have figured out 
how to calculate $\sigma_{\lambda} \ad \f_{\mu}$ for ordered triples.

\ignore {

One more thing I'll need is due to Tits, but formulated more explicitly
in Chapter 4 (Theorem 4.1.2 (ii)) of \refto{Carter}{1972} and as Lemma 2.5 in \refto{Casselman}{2015}:

theorem{tits-geometric-lemma}{If $(\lambda,\mu,\nu)$
is a Tits triple then
$$
	{ p_{\lambda,\mu} + 1 \over \Vert \nu \Vert^{2} } 
	=
	{ p_{\mu,\nu} + 1 \over \Vert \lambda \Vert^{2} } 
	=
	{ p_{\nu,\lambda} + 1 \over \Vert \mu \Vert^{2} } 
$$
}

}

In other words, $p_{\lambda,\mu}$ satisfies a twisted cyclic symmetry.

\section{Computation I}
How do results in the previous section apply
to practical computation of structure constants?

The ultimate goal is to come up with a procedure to determine brackets
$[\f_{\lambda},\f_{\mu}]$ easily, given an invariant basis
$(\f_{\lambda})$.  There are three possibilities.
(1) If $\lambda = -\mu$, the bracket is $h_{\mu}$. 
We can express it as a linear combination of basis elements $h_{\alpha}$:
$$ h_{\mu} = \sum_{\alpha} c_{\alpha}  h_{\alpha} \, , $$
in which the coefficients $c_{\alpha}$ are found in the course
of constructing the roots, since this equation is equivalent to
$$ \mu^{\Vee} = \sum_{\alpha} c_{\alpha} \alpha^{\Vee} \, . $$
(2) The sum $\lambda+\mu$ is not a root,
and the bracket is $0$.
We can be decided by a look-up table of roots.
(3) We have an equation
$$ [\f_{\lambda},\f_{\mu}] = N_{\lambda,\mu} \f_{\lambda+\mu} $$
for some constant $N_{\lambda,\mu}$ of the form $\pm (p_{\lambda,\mu}+1)$.
So we would be given a Tits triple $(\lambda,\mu,\nu)$.
We can rotate it to make it an ordered triple.  According
to \lreftosatz{abc-cor}, our problem
is thus reduced to finding just the values $N_{\lambda,\mu}$ 
when $(\lambda, \mu, \nu)$
is an ordered triple.  Because of invariance under $\theta$,
 $N_{-\lambda,-\mu} = N_{\lambda,\mu}$, and
we may restrict to the case $\lambda > 0$.  

We can in fact
calculate and then store all such values.  
The amount of storage required is roughly proportional
to the number of Tits triples.  As reported in
\refto{Cohen-Murray-Taylor}{2005},
this is of order $r^{3}$,
where $r$ is the rank of the system, so this procedure is entirely
feasible, and noticeably better in storage use than
storing all the $N_{\lambda,\mu}$, since there are roughly $r^{4}$
such pairs.  (Of course using the smaller table
involves more computation.  The trade-off of time versus memory that we see here
is a basic problem in all programming.)

There are three steps to this computation.

\stepno=0
\step In the first, we construct the root system,
without reference to a Lie algebra.
This includes (i) root lengths $\Vert \lambda \Vert$, 
(ii) values of $\< \lambda,\alpha^{\Vee}>$, 
(iii) root reflection tables $s_{\alpha}\lambda$, 
(iv) an expression for each root as a linear combination of the $\alpha$
in $\Delta$, and (v) a corresponding expression for each $\lambda^{\Vee}$
as a sum of $\alpha^{\Vee}$.  
We can also construct a table recording whether or not a given 
array of coordinates is that of a root or not.

\step  
In some way specified in the next sections, we then find
an invariant basis $(\f_{\lambda})$.
It is here where Kottwitz' contribution appears. It will give us also
the associated Tits section $\Dotw$ from $W$ to $N_{G}(T)$,
in which $\Dots_{\alpha}$ for $\alpha$ in $\Delta$.
It will give us at the same 
time all the constants $c(s_{\alpha}, \lambda)$
(with $\alpha$ simple) such that
$$ \Dots_{\alpha} \k{-1} \ad \f_{\mu} = c(s_{\alpha}, \lambda) \f_{s_{\alpha}\mu} \, . $$
Here $\Dots_{\alpha}$ is the same as the `reflection' $\sigma_{\alpha}$
determined by $\f_{\alpha}$ according to Tits'
compatibility.  Miraculously:

\ignore{

{

\parindent=2em

\narrower
\noindent  {\sl Constructing 
the invariant basis $(\f_{\lambda})$ will give at the same time 
formulas for the constants $c(s_{\alpha}, \lambda)$
(with $\alpha$ simple) such that }
$$ \Dots_{\alpha} \k{-1} \ad \f_{\mu} = c(s_{\alpha}, \lambda) \f_{s_{\alpha}\lambda} \, . $$

}

}

I repeat: we start with a frame $(e_{\alpha})$,
but the new basis elements $\f_{\alpha})$
will be different, and the elements $\dots_{\alpha}$
will be different from the $\Dots_{\alpha}$.

We now have a simple recipe for
computing any $\Dotw \k{1} \ad \k{1} \f_{\lambda}$, since if
$$ \Dotw \k{1} \ad f_{\lambda} = c(w, \lambda) \f_{w\lambda} $$
then
$$ c(xy, \lambda) = c(x, y\lambda) c(y,\lambda) \, . $$

{\bold Remark.}  
This can be somewhat inefficient, since the element $w$
can have length up to the number of positive roots.
There is a possible improvement, however, offering a trade of 
memory for time.  Choose an ordering of $\Delta$, and let $W_{i}$ be the subgroup of
$W$ generated by the $s_{\alpha_{j}}$ for $j \le i$.  
As Fokko du Cloux pointed out,
every $w$ in $W$ can be expressed as a unique product
$$ w = w_{1}w_{2} \ldots w_{r} $$
with each $w_{i}$ a distinguished representative of $W_{i-1}\\W_{i}$.
The sizes of these cosets are relatively small, and it is perhaps not infeasible
to store values of the $w\lambda$ and the $c(w, \lambda)$ for $w$ 
a distinguished element in one of them. 

\endremarks

\step
Given the results of the previous step,
we want now to tell how to compute the constants $N_{\lambda, \mu}$
when $(\lambda,\mu,\nu)$ make up an ordered Tits triple
with $\lambda > 0$.

We can do this by
a kind of induction on $\lambda$.  Every positive root $\lambda = w\alpha$
for $w$ in $W$ and $\alpha$ simple.  The \defn{depth} $n$ of $\lambda$
is the minimal length of a chain
$$ \alpha = \lambda_{0} - \lambda_{1} - \cdots - \lambda_{n} = \lambda $$
in which each $\lambda_{i+1} = s_{\alpha_{i}} \lambda_{i}$
for some simple $\alpha_{i}$.  Finding such chains
for all positive roots is part of the natural 
process for constructing the set of roots
in the first place.
If
$$ [ \f_{\lambda}, \f_{\mu}] = N_{\lambda, \mu} \k{1} \f_{-\nu} $$
then
$$ [ \Dots_{\alpha} \ad \f_{\lambda}, \Dots_{\alpha} \ad \f_{\mu}] 
	= N_{\lambda, \mu} ( \Dots_{\alpha} \ad \f_{-\nu} ) \, , $$
and hence
$$ N_{s_{\alpha}\lambda, s_{\alpha}\mu} 
	= c(s_{\alpha},\lambda) c(s_{\alpha},\mu) c(s_{\alpha},-\nu) N_{\lambda, \mu} \, . $$

Reflections transform ordered triples to ordered triples.
Hence if we know how to deal with the case in which $\lambda=\alpha$ is simple
we can compute all the constants for ordered triples in which $\lambda > 0$
by following up the chain.
Furthermore, according to \lreftosatz{tits-tfae} it is very easy to list
ordered triples $(\alpha, \mu, \nu)$.  

Now according to \lreftosatz{rotation} we have, for an ordered triple,
$$ [ \f_{\alpha}, \f_{\mu}] = c(s_{\alpha}, \mu) (-1)^{p_{\alpha,\mu}} 
	(p_{\alpha,\mu}+1) \f_{-\nu} \, . $$
Since $\<\mu, \alpha^{\Vee}> = -1$
we know that $p_{\alpha,\mu}$ is $0$ if $\mu - \alpha$ is not a root,
and is $1$ otherwise (in which case we are dealing with $G_{2}$).

At the end we have the structure constants for all
ordered Tits triples with $\lambda$ positive.  I summarize:

\theorem{summary1}{Suppose that $(\f_{\lambda})$ is an invariant basis,
$\Dotw$ the corresponding Tits section.  Suppose that 
$$ \Dots \ad \f_{\mu} = c(s_{\alpha}, \mu) \f_{s_{\alpha}\mu} $$
for every $\alpha$ in $\Delta$, root $\mu$.  
Suppose $(\lambda, \mu, \nu)$ to be an ordered Tits triple
with $\lambda > 0$.

If $\lambda = \alpha$  lies in $\Delta$, then
$$ N_{\alpha, \mu} = c(s_{\alpha}, \mu) (-1)^{p_{\alpha,\mu}} 
	(p_{\alpha,\mu}+1) \, . $$

Otherwise, we can find $N$ by the induction formula
$$ N_{s_{\alpha}\lambda, s_{\alpha}\mu} 
	= c(s_{\alpha},\lambda) c(s_{\alpha},\mu) c(s_{\alpha},-\nu) N_{\lambda, \mu} 
\, . $$
}

\def\Vee{{\scriptscriptstyle \vee}}
\def\vare{\varepsilon}

\section{Kottwitz' splittings}
{\sl 
It remains to explain how to
construct an invariant basis $(\f_{\lambda})$ and 
give formulas for the constants $c(s_{\alpha}, \lambda)$
(with $\alpha$ simple) such that
$$ \Dots_{\alpha} \k{-1} \ad \f_{\mu} = c(s_{\alpha}, \lambda) \f_{\lambda} 
\quad (\lambda = s_{\alpha}\mu) \, . $$
Here $\Dots_{\alpha}$ is the element of $M_{\alpha}({\Bbb Z})$
compatible with $\f_{\alpha}$.}

In any method of computation in Lie algebras,
the first---and perhaps most important---step is to 
specify an integral basis of the algebra.
\refto{Cohen-Murray-Taylor}{2005} specifies such a basis
in terms of an ordered decomposition of a given root as a sum of
simple ones.  First of all, they
assign an order to the simple roots.   Every positive
root may be expressed uniquely as $\mu = \alpha + \lambda$
in which the height of $\lambda$ is less than that of $\mu$,
and $\alpha$ is least with this property.  They then define
the elements $e_{\mu}$ by induction:
$$ [ e_{\alpha}, e_{\lambda}] = (p_{\alpha,\lambda} + 1) e_{\mu} \, . $$
In effect, such a basis is determined by a choice of spanning
tree in a graph whose nodes are the positive roots,
with a link between each pair $\lambda$ and $\alpha + \lambda$.

The method
I described in \refto{Casselman}{2015a} and \refto{Casselman}{2015b}
chooses a basis in terms of paths in a spanning tree in a different
graph whose nodes are again the positive roots. 
The simplest implementation starts also with an ordering of simple roots.
Every positive root may be expressed as $\mu = s_{\alpha} \lambda$,
with $\lambda$ of smaller height and $\alpha$ minimal.  Then define by induction
$$ e_{\mu} = \dots_{\alpha} e_{\lambda} \, . $$

There is a great deal of arbitrariness in both methods,
since they depend on a somewhat arbitrary choice of spanning tree in a graph.
Kottwitz' contribution is to remove nearly all this annoying
ambiguity.  A basis chosen directly by his method will not be
invariant under $\theta$, but it will be easy to determine from it
one that is.

The original choice of frame gives us Tits' map
$w \mapsto \dotw$ from $W$ back to 
${\cal N}({\Bbb Z})$, and then to the extended group ${\cal N}_{\rm ext}({\Bbb Z})$.
How can it be modified to become a homomorphism?

We are looking for a splitting of the sequence
$$ 1 \longrightarrow \Ss({\Bbb Z}) 
\longrightarrow {\cal N}_{\rm ext}({\Bbb Z})
\longrightarrow W
\longrightarrow 1 \, . $$
This will be of the form
$$ w \longmapsto \Ddotw = \dotw \Cdot \tau_{w} \, , $$
with each $\tau_{w}$ in $\Ss({\Bbb Z})$.
Thus for each root $\beta$ we are looking for a factor
$\tau_{w}(\beta) = \pm 1$.  The map $w \mapsto \Ddotw$
will be a homomorphism if and only if (for $\alpha$ in $\Delta$)
$$ \eqalignno {
	(a)\phantom{'} \k{10} & \Ddotl = 1 \cr
	(b)\phantom{'} \k{10} & \Ddots_{\alpha} \Ddotx 
		= (s_{\alpha}x)\k{1}\Ddotnull \; \hbox{ if $s_{\alpha}x > x$ } \cr
	(c)\phantom{'} \k{10} & \Ddots_{\alpha}\Ddots_{\alpha} = 1 \, . \cr
%
\noalign{These translate directly to properties of $\tau_{w}$:}
%
	(a') \k{10} & \tau_{1} = 1 \cr
	(b') \k{10} & \tau_{s_{\alpha}}(y\beta) \tau_{y}(\beta) 
			= \tau_{s_{\alpha}y}(\beta) 
		\; \hbox{ for all } \beta \; \hbox{ if } \; s_{\alpha}y > y \cr
	(c') \k{10} & (-1)^{\<\beta, \alpha^{\vee}>} 
		= \tau_{s_{\alpha}}(s_{\alpha}\beta) \Cdot \tau_{s_{\alpha}}(\beta) \, . \cr
} $$

We shall see a bit later a fourth 
useful condition on $\Ddotw$ and hence also on $\tau_{w}$.

At any rate, here is Kottwitz' solution of the problem.  For $w$ in $W$ set
$$ 
	R_{w} = \{ \lambda > 0 \, | \, w\lambda < 0 \} \, . 
$$
Thus $\ell(xy) = \ell(x) + \ell(y)$ if and only if
$$ 
	R_{xy} = R_{y} \sqcup y^{-1}R_{x} \, , \equation{rsw} 
$$
and in particular
$$ \eqalign {
	R_{1} &= \emptyset \cr
	R_{s_{\alpha}} &= \{ \alpha \} \cr
	R_{s_{\alpha}w} &= R_{w} \sqcup \{ w^{-1}\alpha \} \quad (w^{-1}\alpha > 0) \, . \cr
} $$

According to Kottwitz' recipe, we have
$$ \tau_{w}(\beta) = (-1)^{F(w,\beta)} \quad \hbox{ with } \quad F(w,\beta) = \sum_{\gamma \in R_{w}} \term{ \beta, \gamma } \, . \equation{ca} $$
The summands are yet to be specified, and
everything in this formula is to be taken modulo $2$.

$\bullet$\ 
Since $R_{1} = \emptyset$ and an empty sum is $0$, 
condition (a) above is immediate.

$\bullet$\ 
What about condition (b)? Suppose $x = s_{\alpha}y > y$.
It must be shown that the cocycle condition
$$ F(s_{\alpha}y, \beta) = F(s_{\alpha}, y\beta) + F(y, \beta) $$
holds. First of all, note that
$$ F(s_{\alpha}, \beta) = \term{\beta, \alpha} $$
since $R_{s_{\alpha}} = \{ \alpha \}$.
Also
$$ F(x, \beta) 
	= \sum_{\gamma \in R_{x}} \term{\beta, \gamma} 
	= \term{\beta, y^{-1}\alpha} + \sum_{\gamma \in R_{y}} \term{\beta, \gamma}
$$
whereas
$$ F(s_{\alpha}, y\beta) + F(x, \beta)\
	= \term{y\beta, \alpha} + \sum_{\gamma \in R_{y}} \term{\beta, \gamma} \, . $$
Therefore (b) will be satisfied if $W$-invariance holds:
$$ \term{w\beta, w\gamma} = \term{\beta, \gamma} \; \hbox{for all $w$ in $W$.} $$

$\bullet$\ 
Condition (c)?  We have
$$ \Ddots_{\alpha} \ad e_{\beta} = (-1)^{\term{\beta, \alpha}} \dots_{\alpha} \ad e_{\beta} \, . $$
Since $\dots^{2}_{\alpha} = \alpha^{\Vee}(-1)$ we thus require that
$$ \term{s_{\alpha}\beta, \alpha} + \term{\beta, \alpha}  
	= \<\beta, \alpha^{\Vee}>\, . $$

This last condition suggests what comes now.
If $\<\beta, \alpha^{\Vee}> = 0$ and hence $s_{\alpha}\beta = \beta$ 
this imposes no condition (since everything is modulo $2$).  
Otherwise $\<\beta, \alpha^{\Vee}>$ 
and $\<s_{\alpha}\beta, \alpha^{\Vee}>$
will be of different signs.  It is therefore natural to set
$$ \term{\beta, \gamma} = \cases {  \<\beta, \gamma^{\Vee}> 
	& if $\<\beta, \gamma^{\Vee}> > 0$ \cr
	\k{10} 0 & if $\<\beta, \gamma^{\Vee}> < 0$.
} \equation{cb} $$
One good sign:

\lemma{w-invariance}{The function $\term{\beta, \gamma}$ is Weyl-invariant.}

\proof/.  Since $\< \beta, \gamma^{\Vee}>$ and $p_{\beta,\gamma}$ are 
both $W$-invariant.\endproof

The requirement that $w\mapsto \Ddotw$ be a homomorphism
imposes no extra condition in the case that $\< \beta, \gamma^{\Vee}> = 0$,
but one more requirement will do so.
I ask now, for reasons that will become apparent in a moment, that 
$$ \Ddotw \ad e_{\beta} = e_{\beta} $$
if $w\beta = \beta$. To guarantee that this occurs, it suffices to 
assume that $\beta$ lies in the closed positive Weyl chamber.
Then the $w$ fixing $\beta$ are generated by simple
root reflections, so we need to require only
that $\Ddots_{\alpha} v_{\beta} = v_{\beta}$ 
($v_{\beta} \in {\goth g}_{\beta}$) for simple roots $\alpha$
with $\<\beta, \alpha^{\Vee}> = 0$. Consideration of
the representation of $\SL_{2}$ corresponding to the root string
tells us that
$$ \dots_{\alpha} \ad  e_{\beta} = (-1)^{p_{\alpha,\beta}} e_{\beta} \, . $$
Therefore
$$ \Ddots_{\alpha}  \ad e_{\beta} 
	= (-1)^{\term{\beta, \alpha}} (-1)^{p_{\alpha,\beta}} e_{\beta} $$
and so we set
$$ \term{\beta, \gamma} = p_{\gamma, \beta} 
\quad \hbox{ if } \< \beta, \gamma^{\Vee}> = 0 \, . \equation {cc} $$
Equations \lreftosatz{cb} and \lreftosatz{cc}
define the terms $\term{\beta, \gamma}$ completely. 
In summary:

\satz{kottwitz}{Theorem}{{\rm (Kottwitz)} Let
$$ \eqalign {
	\term{\beta, \gamma} &= \cases { \< \beta, \gamma^{\Vee}> & if this is positive 
			\raise-4pt\hbox{$\mathstrut$}\cr
	\k{4} p_{\gamma, \beta} & if $\< \beta, \gamma^{\Vee} > = 0$ \cr
	\k{4} \phantom{\Big|} 0 & otherwise. \cr
} \cr
 F(w, \beta) &= \sum_{\gamma \in R_{w}} \term {\beta, \gamma} \cr
	\tau_{w} &= \big( (-1)^{F(w,\beta} \big)_{\beta \in \Sigma} \, . \cr
} $$
Then 
$$ \Ddotw = \dotw \Cdot \tau_{w} $$
is a splitting homomorphism of ${\cal N}_{\ext}({\Bbb Z})$.
In addition, if $w \gamma = \gamma$ then 
$\Ad(\Ddotw)$ is the identity on ${\goth g}_{\gamma}$.}

If the root system is simply laced 
or equal to $G_{2}$ then $s_{\lambda}\beta = \beta$
implies that $p_{\lambda,\beta} = 0$.  Therefore the non-trivial case occurs only
for systems $B_{n}$, $C_{n}$, or $F_{4}$.

{\bold Remark.} {\sl 
Lemma 2.1A of \refto{Langlands-Shelstad}{1987} exhibits the $2$-cocycle defining
the extension ${\cal N}({\Bbb Z})$ 
determined by Tits's splitting $w \mapsto \dotw$.  Explicitly,
$$ \dotx \doty = \kappa(x, y) (xy)\k{2}\dott \; \hbox{ with } \; \kappa(x,y) 
	= \k{-10} \prod_{{\scriptstyle \gamma > 0 \atop \scriptstyle x^{-1}\gamma < 0 } \atop
	\scriptstyle y^{-1}x^{-1}\gamma < 0} \k{-10} \gamma^{\Vee}(-1) \, . $$
Does Kottwitz' splitting allow arguments 
of Langlands and Shelstad to be simpler?}

\endremarks

The $W$-orbits in $\Sigma$ are the sets of all roots of the same length.  
Pick one simple root $\alpha$ in each orbit, and let $\e_{\alpha} = e_{\alpha}$
be the corresponding element in the frame chosen at the beginning.
If $\lambda = w\alpha$ is root 
with $\alpha$ equal to one of these distinguished choices, define 
$$ \e_{\lambda} = \Ddotw  \ad \e_{\alpha} \, . $$
The definition of $F(s_{\alpha}, \beta)$
in the case when $\<\beta, \alpha^{\Vee}> = 0$ insures that this
is a valid definition.  As a consequence of \lreftosatz{kottwitz}:

\satz{kottwitz-basis}{Corollary}{The 
integral basis $(\e_{\gamma})$ of 
$V_{\Bbb Z}$ is such that $\Ddotw  \ad \e_{\gamma} = \e_{w\gamma}$
for all roots $\gamma$ and $w$ in $W$.}

\definition{semi-can}{I'll call such a basis \defn{semi-canonical}.}

There are a small set of possibilities, two for each $W$-orbit
in $\Sigma$. 

\ignore {

In practice, we shall want to compute $\taveu_{w}$ 
explicitly only when $w = s_{\alpha}$ for $\alpha$ in $\Delta$.  In this case,
there is a simplification, since $R_{s_{\alpha}}$ is a singleton.

proposition{singleton}{For $\alpha$ in $\Delta$,
$x_{\lambda}$ in ${\goth g}_{\lambda}$
$$ \Ddots_{\alpha} \ad x_{\lambda} = (-1)^{\term{\lambda, \alpha}} x_{\lambda} \, . $$
}

}

{\bold Example.} For a simply laced root system, if $\< \beta, \alpha^{\Vee}> = 0$
then $p_{\alpha,\beta} = 0$.  Therefore
$$ \tau_{s_{\alpha}}(\lambda) = \cases {
	(-1)^{\< \beta, \alpha^{\Vee}>} & if $\< \lambda, \alpha^{\Vee}> > 0$ \cr
		\k{12}1 & otherwise. \cr
} $$
This applies in particular to $G = \SL_{3}$.
Take $\alpha$, $\beta$ as the standard simple roots,
and let $\gamma = \alpha+\beta$.  
Recall that $e_{i,j}$ is the matrix with a single 
non-zero entry $1$ at 
$(i,j)$.  Choose $e_{1,2}$ and $e_{2,3}$ to define the frame,
spanning the root spaces for $\alpha$, $\beta$. 
The corresponding elements of ${\cal N}({\Bbb Z})$ are
$$ 
	\dots_{\alpha} = \mmatrix { \k{-3}\phantom{-}\zero & 1 & \zero \cr
					 \k{-3}-1 & \zero & \zero \cr
					 \k{-3}\phantom{-}\zero & \zero & 1 \cr }, \quad 
	\dots_{\beta} = \mmatrix { 1 &  \k{-4}\phantom{-}\zero & \zero \cr
					\zero &  \k{-4}\phantom{-}\zero & 1 \cr
					\zero &  \k{-4}-1 & \zero \cr } \, . 
$$
And here is a table of the $\ad$ actions:
\medbreak

\settabs
\+ xxxxxxxxx 
& xxxx   
& xxxxxx   
& xxxxxxxxxxxx 
& xxxxxxxx 
& xxxxxxxx 
& xxxxxxxxxxxx 
& xxxxxxxx 
& xxxxxxxx 
& \cr

\+ 
&  $\phantom{-}\lambda$ 
& \hfill $e_{\lambda\phantom{.a}}$ 
& \hfill $\dots_{\alpha} \ad e_{\lambda}$\k{-4}  
& \hfill $\<\lambda,\alpha^{\Vee}> $\k{-6} 
& \hfill $\term{\lambda, \alpha}$\k{-6}
& \hfill $\dots_{\beta} \ad e_{\lambda}$\k{-4}  
& \hfill $\<\lambda,\beta^{\Vee}>$\k{-6} 
& \hfill $\term{\lambda, \beta}$\k{-6}
& \cr
\+ & $\phantom{-}\alpha$ 
& \hfill $e_{1,2}$ 
& \hfill $-e_{2,1}$ 
& \hfill $\ph 2$ 
& \hfill $0$ 
& \hfill $-e_{1,3}$ 
& \hfill $-1$ 
& \hfill $0$ 
& \cr
\+ 
& $\phantom{-}\beta$ 
& \hfill $e_{2,3}$ 
& \hfill $\ph e_{1,3}$ 
& \hfill $-1$ 
& \hfill $0$ 
& \hfill $-e_{3,2}$ 
& \hfill $\ph 2$ 
& \hfill $0$ 
& \cr
\+ 
& $\phantom{-}\gamma$ 
& \hfill $e_{1,3}$ 
& \hfill $- e_{2,3}$ 
& \hfill $\ph 1$ 
& \hfill $1$ 
& \hfill $\ph e_{1,2}$ 
& \hfill $\ph 1$ 
& \hfill $1$ 
& \cr
\+ 
& $-\alpha$ 
& \hfill $e_{2,1}$ 
& \hfill $-e_{1,2}$ 
& \hfill $-2$ 
& \hfill $0$ 
& \hfill $- e_{3,1}$ 
& \hfill $\ph 1$ 
& \hfill $1$ 
& \cr
\+ 
& $-\beta$ 
& \hfill $e_{3,2}$ 
& \hfill $\ph e_{3,1}$ 
& \hfill $\ph 1$ 
& \hfill $1$ 
& \hfill $-e_{2,3}$ 
& \hfill $-2$ 
& \hfill $0$ 
& \cr
\+ 
& $-\gamma$ 
& \hfill $e_{3,1}$ 
& \hfill $-e_{3,2}$ 
& \hfill $-1$ 
& \hfill $0$ 
& \hfill $\ph e_{2,1}$ 
& \hfill $-1$ 
& \hfill $0$ 
& \cr
\medbreak
 
If we start with $\e_{\alpha} = e_{1,2}$ we get
$$ \eqalign {
	\e_{\alpha} &= e_{1,2} = e_{\alpha} \cr
	\e_{\gamma} &= \Ddots_{\beta} \e_{\alpha} \cr
				&= (-1)^{0} \dots_{\beta} \ad e_{1,2}  \cr
				&= -e_{1,3} \cr
	\e_{\beta} &= \Ddots_{\alpha} \e_{\gamma} \cr
				&= (-1)^{1} \dots_{\alpha} \ad (-e_{1,3}) \cr
				&= -e_{2,3} = -e_{\beta} \, . \cr
} $$

Thus:

\proposition{sl3-e}{If\/ $G = \SL_{3}$ and 
$\e_{\alpha} = e_{\alpha}$, then $\e_{\beta} = -e_{\beta}$.}

This example has consequences for arbitrary root systems.

Something very similar is true for all groups $\SL_{n}$.  Here, choose
the base point of the Dynkin diagram to be the end point
corresponding to the simple root $\varepsilon_{1} - \varepsilon_{2}$.
Then $\e_{i,j} = (-1)^{j} e_{i,j}$.

\endremarks

A semi-canonical basis will not be invariant under $\theta$,
but it is easy to see how it fails, and then how to modify it to be so.
Recall that the height of a root is defined by the formula
$$ \height\Big( {\sum}_{\Delta} \lambda_{\alpha} \alpha \Big) =
	{\sum}_{\Delta} \lambda_{\alpha} \, . $$

\theorem{ht-parity}{For any root $\gamma$ and Kottwitz basis $(\e_{\gamma}$) 
$$ \e^{\theta}_{\gamma} = (-1)^{\height(\gamma)-1} \e_{-\gamma} \, . $$
In particular, if $\alpha$ is simple then
$$ \e_{\alpha}^{\theta} = \e_{-\alpha} \, . $$
}

This particularly simple formulation is due to Kottwitz.

\proof/.  In a number of short steps.

\stepno=0
\step
The following is straightforward:

$\bullet$\ {For all $\beta$, $\gamma$
$$ \dbll \beta, \gamma \dblr + \dbll -\beta, \gamma \dblr 
		= \< \beta, \gamma^{\Vee}> $$
}

This is to be interpreted modulo $2$, of course.  

\step
Now let
$$ h(w, \beta) = \sum_{\gamma \in R_{w}} \< \beta, \gamma^{\Vee}> \, . $$

$\bullet$\ {\sl For $v$ in ${\goth g}_{\beta}$
$$ (\Ddotw  \ad  v )^{\theta} = (-1)^{h(w,\beta)} \, \Ddotw \ad v^{\theta}  \, . \equation{theta-Ddot} $$
}

This is because $\dots_{\alpha}^{\theta} = \dots_{\alpha}$.

\step
Induction on the length of $w$ together with \lreftosatz{rsw} will prove:

$\bullet$\ {\sl For $w$ in $W$ and root $\lambda$
$$ \height(w\lambda) - \height(\lambda) 
	= h(w, \lambda) \, . $$
}

This concludes the proof of the Theorem.\endproof

In order to specify the $\e_{\gamma}$, given a frame $(e_{\alpha})$,
we fix one simple root $\alpha$ in each $W$-orbit,
and set $\e_{\alpha} = e_{\alpha}$.  Fixing the $\e_{\beta}$
for other simple roots $\beta$ is then very easy.
For finite-dimensional Lie algebras, 
$W$-orbits of roots are in correspondence with possible
root lengths.  For irreducible
systems, there are at most two possible lengths,
and the simple roots of a given length make up
a connected segment $\Xi$ in the Dynkin diagram. It is only in systems
$B$, $C$, $F$, and $G$ that there are two lengths,
and only for system $F$ is there more than one simple
root of each length.

To determine the elements $\e_{\lambda}$ choose,
somewhat arbitrarily, one \defn{special} root  $\alpha_{\Xi}$ 
on each segment $\Xi$.  
For every simple root $\alpha$, let 
$$ d(\alpha) = \hbox{ the distance from $\alpha$ 
					to the special root $\alpha_{\Xi}$ in its segment.} $$
Any two neighbours in the Dynkin diagram of
the same length lie in the simple root system of a copy of $\SL_{3}$.
The choice of $\e_{\alpha}$ determines an element $\sigma_{\alpha}$.  
The following is a consequence of \lreftosatz{sl3-e}:

\corollary{ea-lemma}{For $\alpha$ in $\Delta$
let $c_{\alpha} = (-1)^{d(\alpha)}$.  Then
$$ \eqalign {
	\e_{\alpha} &= c_{\alpha} \, e_{\alpha} \cr
	\sigma_{\alpha} &= \alpha^{\Vee}(c_{\alpha}) \, \dots_{\alpha} \, . \cr
} $$
}

Here, I recall, $\sigma_{\alpha}$ is the element of $N_{\Bbb Z}(T)$
associated by Tits' scheme to the choice of $\e_{\alpha}$
as basis of ${\goth g}_{\alpha}$ (or of $\e_{-\alpha}$ for ${\goth g}_{-\alpha}$).

{\bold Remark.} \lreftosatz{sl3-e} is the obstruction to
extending these the results to Kac-Moody algebras
whose Dynkin diagram has loops containing an odd number
of roots connected by simple edges.

{\bold Remark.}  I have mentioned the 'root graph'
without being precise, and I should say something more about it.
It is a graph whose nodes are the positive roots, and its base is made up 
of the simple roots.  There is an oriented edge from $\lambda$ to $s_{\alpha}\lambda$
if and only if $s_{\alpha}\lambda$ has greater height than
$\lambda$, or equivalently if and only if $\< \lambda,\alpha^{\Vee}> < 0$.
This is very useful, since in these circumstances $\term{\lambda,\alpha}$
is always $0$.  One consequence is an easy construction of
the basis $(\e_{\lambda})$.  Following upward
links in the root graph, one represents every root as
 an increasing chain
$$ \alpha = \lambda_{0} \prec \ldots \prec \lambda_{n} = \lambda 
\quad (\lambda_{i+1} = s_{\alpha_{i}} \lambda_{i}) $$
and then 
$$ \e_{\lambda} = \Ddots_{n-1} \ldots \Ddots_{0} \ad \e_{\alpha} \, . $$
This is very useful for debugging programs, since for the classical root
systems one can construct Kottwitz' basis in terms of explicit matrices,
for which one can calculate Lie brackets in terms of matrix products.

\section{Computation II}
Define
$$ \gamma(\lambda) =\cases {
		\k{11} 1 \k{-11} \phantom{(-1)^{\height(-\lambda)-1}} & if $\lambda > 0$ \cr
		(-1)^{\height(-\gamma)-1} & if $\lambda < 0$. \cr
}  $$
As an immediate consequence of \lreftosatz{ht-parity}:

\proposition{T-basis}{Given the Kottwitz basis $(\e_{\lambda})$,
the elements 
$$ \f_{\lambda} = \gamma(\lambda) \e_{\lambda} $$
form an invariant integral basis.} 

{\bold Remark.}  I emphasize: we start
with a given frame, then find a new frame that is rarely the same
as the original.  It is this new frame that we extend to an integral basis in
a uniquely determined way.

\endremarks

Let $\Dotw$ be the corresponding Tits section. 

I want to summarize some of our current situation. We started with
a frame $(e_{\alpha})$, with associated elements $\dots_{\alpha}$.
We then defined Kottwitz' section $\Ddotw$, and from this
an integral basis $\f_{\lambda}$.  To the $\f_{\alpha}$ we then have 
new Tits sections $\Dots_{\alpha}$.  This gives us
$$ \eqalign {
	\Ddots_{\alpha} \ad x_{\lambda} 
		&= (-1)^{\term{\lambda, \alpha}} \dots_{\alpha} \ad x_{\lambda} \cr
	\Dots_{\alpha} \ad x_{\lambda} 
		&= (c_{\alpha})^{\< \lambda, \alpha^{\Vee}>} \dots_{\alpha} \ad x_{\lambda} \, . \cr
} $$
I recall that
$$ c_{\alpha} = (-1)^{d(\alpha)} \, , $$
where $d(\alpha)$ for $\alpha$ in $\Xi$ 
measures distance along the Dynkin diagram
from the nearest simple root $\alpha_{\Xi}$.


The following result
encapsulates the basic reason why Kottwitz' basis makes computation simple.

\satz{basic-formula}{Theorem}{Suppose $\alpha$ to lie in $\Delta$.  If
$$ m_{\alpha,\lambda} = 
		(-1)^{\term{\lambda, \alpha}} c_{\alpha}^{\<\lambda,\alpha^{\Vee}>} $$
for $\lambda > 0$, then for every pair $\lambda > 0$ and $\alpha$ in $\Delta$ 
$$ \Dots_{\alpha} \ad \f_{\lambda} = c(s_{\alpha}, \lambda) \f_{s_{\alpha}\lambda} $$
with
$$
c(s_{\alpha}, \lambda)= \cases {
	m_{\alpha,\lambda} \k{4}  & if $\lambda > 0$ \cr
	m_{\alpha,-\lambda}  	& if $\lambda < 0$ . \cr
} $$
}

\proof/.  Since $\lambda(\alpha^\Vee(x)) = x^{\< \lambda, \alpha^{\Vee}>}$:
$$ \eqalignno {
	\Ddots_{\alpha} \ad \e_{\lambda}
		&= \e_{\mu} \cr
		&= (-1)^{\term{ \lambda, \alpha}} \, \dots_{\alpha} \ad \e_{\lambda} &  \cr
	\dots_{\alpha} \ad \e_{\lambda} 
		&= c_{\alpha}^{\<\lambda,\alpha^{\Vee}>} \Dots_{\alpha} \ad \e_{\lambda} \cr
	\Dots_{\alpha} \ad \e_{\lambda}
		&= c_{\alpha}^{\<\lambda,\alpha^{\Vee}>} (-1)^{\term{ \lambda, \alpha}} \e_{\mu} \cr
		&= m_{\alpha,\lambda} \Cdot \e_{\mu} \, . \cr
} $$
This concludes when $\lambda > 0$, even if $\lambda = \alpha$
and $\sigma_{\alpha} \ad \k{0.8} \f_{\alpha} = \f_{-\alpha}$.  When not, 
apply the involution $\theta$ to this equation, noting that $\Dots_{\alpha}$
commutes with it.\endproof

{\bold Example.} Look at $\SL_{3}$ again.  What is $[\e_{\alpha}, \e_{\beta}]$?
$$ \eqalign {
	c_{\alpha} &= \phantom{-} 1 \cr
	\< \beta, \alpha^{\Vee} > &= -1 \cr
	\term{\alpha,\beta} &= \phantom{-} 0 \cr
	p_{\alpha,\beta} &= \phantom{-} 0 \cr
	c(s_{\alpha}, \beta) &= \phantom{-} 1 
} $$
Hence $\Dots_{\alpha}\e = \e_{\gamma}$.

{\bold Example.} Say $G = \Sp(4)$.  Let $\alpha = \varepsilon_{0}-\varepsilon_{1}$
and $\beta = 2\varepsilon_{1}$ be the simple roots.  
Since there are two lengths of roots, we may set as frame
$$
\e_{\alpha} = e_{\alpha} = \mmatrix {
0 & 1 & 0 & \phantom{-}0 \cr
0 & 0 & 0 & \phantom{-}0 \cr
0 & 0 & 0 & -1 \cr
0 & 0 & 0 & \phantom{-}0 \cr
},\quad 
\e_{\beta} = e_{\beta} = \mmatrix {
0 & 0 & 0 & 0 \cr
0 & 0 & 1 & 0 \cr
0 & 0 & 0 & 0 \cr
0 & 0 & 0 & 0 \cr
} \, . 
$$
Then
$$ \dots_{\beta} =   
\mmatrix {
1 & \phantom{-}0 & 0 & 0 \cr
0 & \phantom{-}0 & 1 & 0 \cr
0 & -1 & 0 & 0 \cr
0 & \phantom{-}0 & 0 & 1 \cr
} \, . 
$$
Since $\<\alpha,\beta^\Vee> = -1$, $s_{\beta}\alpha = \alpha+\beta = \hbox{ (say) } \gamma$.
Also, $\term{\beta,\alpha} = 0$ and hence
$$ \e_{\gamma} = 
\dots_{\beta} \e_{\alpha} \dots_{\beta}^{-1} 
= \mmatrix {
1 & \phantom{-}0 & 0 & 0 \cr
0 & \phantom{-}0 & 1 & 0 \cr
0 & -1 & 0 & 0 \cr
0 & \phantom{-}0 & 0 & 1 \cr
} 
\mmatrix {
0 & 1 & 0 & \phantom{-}0 \cr
0 & 0 & 0 & \phantom{-}0 \cr
0 & 0 & 0 & -1 \cr
0 & 0 & 0 & \phantom{-}0 \cr
}
\mmatrix {
1 & 0 & \phantom{-}0 & 0 \cr
0 & 0 & -1 & 0 \cr
0 & 1 & \phantom{-}0 & 0 \cr
0 & 0 & \phantom{-}0 & 1 \cr
}
=
\mmatrix {
0 & 0 & -1 & \phantom{-}0 \cr
0 & 0 & \phantom{-}0 & -1 \cr
0 & 0 & \phantom{-}0 & \phantom{-}0 \cr
0 & 0 & \phantom{-}0 & \phantom{-}0 \cr
} \, . $$
One calculates directly that
$$ [\e_{\alpha},\e_{\beta}] = \mmatrix {
0 & 0 & 1 & 0 \cr
0 & 0 & 0 & 1 \cr
0 & 0 & 0 & 0 \cr
0 & 0 & 0 & 0 \cr
} = -\e_{\gamma} \, . $$
But it is instructive to trace how the computations in this paper would go.
We are looking at the triple
$(\alpha,\beta, -\gamma)$.  Since $\Vert \beta \Vert = 2$
while $\Vert \alpha \Vert = 1$, the associated ordered triple is
$(\beta, -\gamma, \alpha)$.  Since $-s_{\beta}\gamma = -\alpha$,
we must next compute the constant $\varepsilon$ such that
$$ \dots_{\beta}\f_{-\gamma} = \varepsilon \f_{-\alpha} \, . $$
This is $c(s_{\beta}, -\gamma)$, which according to \lreftosatz{basic-formula} is
$$ m_{\beta, \gamma} = (-1)^{\term{\gamma,\beta}} = -1 \, . $$

{\bold Remark.}  It is not difficult to compute any $p_{\lambda,\mu}$
by finding directly the maximum value of $n$ such that
$\mu - n\lambda$ is a root.  But this is more expensive in time
than necessary.  The circumstances in which we have to compute $p_{\lambda,\mu}$ are 
in fact somewhat limited: (1) when $\lambda = \alpha$
is simple and $\< \mu, \alpha^{\Vee} > = -1$;
(2) when $\alpha$ is simple and $\< \mu, \alpha^{\Vee}> = 0$; (3) $(\lambda, \mu, \nu)$
form a Tits triple.  In case (1) or (2), we just have to check whether
$\mu - \alpha$ is a root.  But if $\alpha$ is simple,
computing $\mu - \alpha$ is trivial, a matter of decrementing one coordinate.
In case (3), we can apply \lreftosatz{tits-geometric-lemma}
in order to reduce to the case in which $(\lambda, \mu, \nu)$ is an ordered
triple.  These are dealt with in the process of ascending the root graph
that is mentioned at the end of \S 3, since
$p_{s_{\alpha}\lambda, s_{\alpha}\mu} = p_{\lambda,\mu}$.

\vfill\eject

\section{Summary}
I summarize here the computation.

We start with a semi-simple Lie algebra ${\goth g}$ over ${\Bbb C}$
and a frame $(e_{\alpha})$.  This is given to us
as a Cartan matrix.  We construct the positive roots.It is best to give for each
its coordinates in two forms for each of $\lambda$ and $\lambda^{\Vee}$,
and also $\Vert \lambda \Vert^{2}$ equal to $1$, $2$, or $3$.

We can calculate the lengths according to the rule
that if $\ell = |\<\alpha, \beta^{\Vee}>| \ge 2$ then $\Vert \alpha \Vert^{2} = \ell$
and $\Vert \beta \Vert^{2} = 1$.

{\sl How to compute the constants $N_{\lambda,\mu}$ with respect to
the invariant basis determined by Kottwitz' method?}

We shall use the formulas given earlier for $\term{\alpha, \mu}$,
$d_{\alpha}$, and $c(s_{\alpha}, \mu)$ when $\alpha$ is a simple root.

To calculate $N_{\lambda, \mu}$, first check whether
$\lambda + \mu$ is a root.  If so, 
compute $P = p_{\lambda, \mu}$ and let $\nu = -(\lambda = \mu)$.
Rotating, we may now assume that $(\lambda, \mu, \nu)$ is
an ordered triple.  If $\lambda < 0$, change signs of all three.

If $\lambda$ is not simple, find a simple root $\beta$ such that
$s_{\beta}\lambda$ has smaller height.  Apply $\dots_{\beta}$ to the triple,
to replace the original triple, which will remain ordered.
Keep track of the factor $c(s_{\beta},\lambda)c(s_{\beta},\mu)c(s_{\beta},\nu)$.
Continue, accumulating products in the factor until
$\lambda$ is simple.  Then apply \lreftosatz{summary1},
replacing $p_{\lambda,\mu}$ by the value $P$.

For speed, it will be useful to store values of
$N_{\lambda, \mu}/(p_{\lambda, \mu}+1)$ whenever $(\lambda, \mu, \nu)$
is an ordered triple.

For debugging, it is useful to check that Jacobi's equation is valid:
$$ [a, [b, c]] + [b, [c, a]] + [ c, [a, b]] = 0 $$
It is also useful to check by comparison with explicit matrix calculations
for classical groups.

\references

\ref/Lisa Carbone, Alexander Conway, Walter Freyn, and Diego Penta,
\article{Weyl group orbits on Kac-Moody root systems},
{\tt arXiv.1407:3375}, 2015.

\ref/Roger W. Carter,
\book{Simple groups of {L}ie type}, 
John Wiley \& Sons, 1972.

\ref/Bill Casselman, 
\article{On Chevalley's formula for structure constants},
\jour{Journal of Lie Theory} \vol{25} (2015),
431--441.

\ref/\again/,
\article{Structure constants of Kac-Moody Lie algebras}, 
in \refto{Howe et al.}{2015}.

\ref/Claude Chevalley,
\article{Sur certains groupes simples},
\jour{T\^ohoku Mathematics Journal} \vol{48} (1955),
14--66.

\ref/Arjeh Cohen,
Scott Murray,
and Don Taylor,
\article{Computing in groups of Lie type},
{\sl Mathematics of Computation} {\bold 73} (2004), 1477--1498.

\ref/
Igor Frenkel and Victor Kac,
\article{Affine Lie algebras and dual resonance models},
\jour{Inventiones Mathematicae} \vol{62} (1980), 23--66.

\ref/Roger Howe \etal/ (editors),
\book{Symmetry: representation theory and its applications}, volume \vol{257}
in the series \jour{Progress in Mathematics}, Elsevier, 2015.

\ref/Robert E. Kottwitz, personal communication in May, 2014.

\ref/R. P. Langlands and D. Shelstad,
\article{On the definition of transfer factors},
\jour{Mathematische Annalen} \vol{278} (1987), 219--271.

\ref/
L. J. Rylands, `Fast calculation of structure constants', preprint, 2000.

\ref/Jacques Tits ([T]),
\article{Sur les constants de structure et le th\'eor\`eme d'existence
des alg\`ebres de Lie semi-simple},
{\sl \frenchspacing Publications de {\frenchspacing l'I. H. E. S.}}
\vol{31} (1966), 21--58.

\ref/\again/,
\article{Normalisateurs de tores I.  Groupes de Coxeter \'etendus},
\jour{Journal of Algebra} \vol{4} (1966), 96--116.

\bye